\documentclass[12pt]{amsart}
\usepackage[cmtip,arrow]{xy}
\usepackage{pb-diagram,pb-xy,
amsxtra,amsmath,amsthm,amssymb,amscd,ascmac,amsfonts,
multicol,tabularx,latexsym,mathrsfs}
\usepackage{amsfonts}
\theoremstyle{definition}

\newcommand{\Q}{\mathbb Q}
\newcommand{\Z}{\mathbb Z}

\newcommand{\Gal}{\mathrm{Gal}}
\newtheorem{thm}{Theorem}

\newtheorem{lem}{Lemma}
\newtheorem{prop}{Proposition}
\newtheorem{rem}{Remark}

\newcommand{\rank}{\mathrm{rank}}

\newcommand{\Cl}{\mathrm{Cl}}
\newdimen\minCDarrowwidth
\date{}
\begin{document}
\title[]
{Construction of a cyclic $p$-extension of number fields whose unit group
has prescribed Galois module structure}
\author{Manabu\ OZAKI}
\maketitle
\section{Introduction}
Let $p$ be a fixed prime number.
For any number field $K$ and finite set $S$ of primes of $K$,
we denote by $E_{K,S}$ the $S$-unit group of 
$K$.
We define $\mathcal{E}_{K,S}:=\left(E_{K,S}/\mathrm{Tor}_{\Z}(E_{K,S})\right)\otimes\Z_p$.
We are interested in the Galois module structure
of $\mathcal{E}_{K,S}$ for Galois extensions $K/k$ of number fields and
Galois stable sets $S$ of primes of $K$,
especially, in what kind of $\Z_p[G]$-modules occur as
$\mathcal{E}_{K,S}$ for $G$-extensions $K/k$
when we fix a finite group $G$.

Burns \cite{B}, Burns-Lim-Maire \cite{BLM}, and 
Kumon-Lim \cite{KL} investigated this problem
and showed that the Galois module structure of $\mathcal{E}_{K,S}$
is highly restricted if we bound the size of $S$,
the number of the ramified primes $R(K/k)$ of $k$ in $K/k$ and
the $p$-rank $r_p(K/k)$ of $S$-ideal class group of 
$K(\mu_p)$ (or the intermediate fields of $K/k$ in \cite{BLM}).

Specifically, they obtained finiteness theorems for the
number of the $\Z_p[G]$-modules without projective direct factors 
which appear as a direct summand of $\mathcal{E}_{K,S}$
when $(K/k, S)$ runs over the pairs of 
$G$-extensions $K/k$ of number fields 
with bounded $\# R(K/k)$ and $r_p(K/k)$, and
the finite sets $S$ of primes of $k$ with
bounded $\#S$.

In the present paper, we will consider the situation
that $G$ is a cyclic $p$-group and $S=\emptyset$, and that 
we make no restrictions on the $G$-extensions $K/k$.
Then, contrary to the above mentioned finiteness theorems,
we will show that every $\Z_p$-torsion free finitely generated
$\Z_p[G]$-module in fact appears as a direct summand of $\mathcal{E}_K:=\mathcal{E}_{K,\emptyset}$:
\begin{thm}
Let $p$ be an odd prime number and 
$G$ a finite cyclic $p$-group. 
Assume that $C$ is a finitely generated $\Z_p[G]$-module
such that $\mathrm{Tor}_{\Z_p}(C)=0$ and
\[
C\otimes_{\Z_p}\Q_p\simeq \Q_p[G]^{\oplus r}\oplus 
I_{\Q_p[G]}
\]
as $\Q_p[G]$-modules for $r\ge 0$,
where $I_{\Q_p[G]}$ stands for the augumentation ideal of
$\Q_p[G]$.
Then there exists an unramified $G$-extension $K/k$
of number fields such that
\[
\mathcal{E}_K\simeq C\oplus\Z_p[G]^{\oplus s}
\]
as $\Z_p[G]$-modules for a certain $s\ge 0$.
\end{thm}   

\

The outline of the present paper is as follows:

In Section 2, we will describe the Galois cohomology groups
of the unit group $E_K$ for unramified extensions $K/k$
in terms of the Galois group of a certain unramified extension
associated to $K/k$.
Then, applying the existence theorem of maximal unramified
$p$-extensions with prescribed Galois groups (\cite{HMR},\cite{Oz}),
we will show that there exists an unramified $G$-extension $K/k$
such that the cohomology groups of $\mathcal{E}_K$
is isomorphic to those of a given $\Z_p[G]$-module $C$
as in the statement of the theorem when $G$ is a cyclic $p$-group.

In Section 3, we will introduce Yakovrev's theory (\cite{Ya})
for the determination of a $\Z_p$-torsion free $\Z_p[G]$-module $M$
modulo free direct factors, $G$ being a cyclic $p$-group,
in terms of a certain diagram consisting of 1-dimensional cohomology
groups of $M$, the restriction and the co-restriction maps.

In Section 4, combining results obtained in Section 2 and 
Yakovrev's theory, we will complete the proof of Theorem 1. 

In Section 5, we will give examples of explicit $\Z_p[\Gal(K/k)]$-module
structures of $\mathcal{E}_K$ for unramified cyclic $p$-extensions
$K/k$ by using the method explored in the preceding sections.
\section{Description of Galois cohomology groups of unit groups.}
%
%
%
%
%
%
%
We recall the notion of a {\it splitting module} to
describe cohomology groups of unit groups in terms of structures
of certain Galois groups.

Let 
\begin{equation*}
(\varepsilon)
\ \ \ \ \ 
1\longrightarrow A\longrightarrow
\mathcal{G}\longrightarrow G\longrightarrow 1
\end{equation*}
be a group extension with abelian kernel $A$.
We denote by $\gamma_{\varepsilon}\in H^2(G,A)$
the cohomology class associated to $(\varepsilon)$.
Define the $\Z[G]$-module $M_{(\varepsilon)}$ so that
\begin{equation*}
M_{(\varepsilon)}
=A\oplus\bigoplus_{1\ne\tau\in G}
\Z b_\tau
\end{equation*}
as $\Z$-modules.
Here $\{b_\tau\}$ is a free $\Z$-basis and
the $G$-action on $M_{(\varepsilon)}$
is given by the natural $\Z[G]$-module structure of $A$
defined by $(\varepsilon)$ and 
\begin{equation*}
\sigma b_\tau=b_{\sigma\tau}-b_{\sigma}
+f(\sigma,\tau)
\end{equation*}
for $\sigma,\,\tau\in G$,
where $f$ is a $2$-cocycle in the cohomology class
$\gamma_{\varepsilon}$
and we set $b_1:=f(1,1)\in A$.
We call $M_{(\varepsilon)}$ {\it the splitting module associated to}
$(\varepsilon)$, and note that 
the $\Z[G]$-module isomorphism class of $M_{(\varepsilon)}$
is independent of the choice of a 2-cocycle $f$
of the cohomology class $\gamma_\varepsilon$.

Then we derive from the group extension $(\varepsilon)$
the exact sequence of $\Z[G]$-modules
\begin{equation*}
(\varepsilon^*)\ \ \ \ \ 0\longrightarrow A\longrightarrow M_{(\varepsilon)}
\overset{\pi}{\longrightarrow} I_{\Z[G]}\longrightarrow 0,
\end{equation*}
where $I_{\Z[G]}$ denotes the augmentation ideal of $\Z[G]$ 
and the map $\pi$
is given by $b_\sigma\mapsto \sigma-1\ \ (\sigma\in G)$.

Now we will give a description of the Galois cohomology group
$\hat{H}^i(\Gal(K/k),E_K)$ for an unramified Galois extension $K/k$ 
in terms of a certain group extension of Galois groups associated to $K/k$. 

Let $K/k$ be an unramified Galois extension of number fields,
and we denote by $H_K$ the maximal unramified abelian extension of $K$.
Put $G:=\Gal(K/k),\ \mathcal{G}:=\Gal(H_K/k)$, and
$A:=\Gal(H_K/K)$.
Then we have the natural group extension 
\[
(\varepsilon_{K/k})\ \ \ \ \ 1\longrightarrow A\longrightarrow
\mathcal{G}\longrightarrow G\longrightarrow 1.
\]
The following proposition gives a description of the Galois cohomology groups of $E_K$ in therms of the group extension $(\varepsilon_{K/k})$:
%
%
%
%
%
%
\begin{prop}\label{prop1}
For any unramified Galois extension $K/k$ of number fields 
with $G=\Gal(K/k)$ and $i\in\Z$, there exist $\Z[G]$-module isomorphisms
\[
\varphi_{H}:\hat{H}^i(H, E_K)\simeq \hat{H}^{i-2}(H,M_{(\varepsilon_{K/k})}).
\]
for the subgroups $H$ of $G$ such that
$\varphi_H$'s are compatible with the restriction and 
the co-restriction maps between cohomology groups.
\end{prop}

{\it Proof.}\ \ \ In fact, this proposition follows from a special case of
the Tate sequence given by Ritter-Weiss \cite{RW1996}.
However we will give a direct proof,
which is based on essentially same method as theirs
specialized to the most simple situation, namely,
the Galois extension considered is unramified. 

Let $J_K$, $U_K$, $C_K$, and $\Cl_K$ be the idele group, the unit idele
group, the idele class group,
and the ideal class group of $K$, respectively.
We will observe the following exact commutative diagram:
\begin{equation*}
\begin{diagram}
\dgARROWLENGTH=0.5cm
\dgHORIZPAD=0.5cm
\node[2]{0}
\arrow{s}
\node{0}
\arrow{s}\\
\node[2]{U_K/E_K}\arrow{e,t}{\sim}\arrow{s}\node{N}\arrow{s}\\
\node{0}\arrow{e}\node{C_K}\arrow{s}\arrow{e}
\node{M_{(\overline{\varepsilon})}}\arrow{s}\arrow{e}\node{I_G}
\arrow{e}\arrow{s,=}\node{0}\\
\node{0}\arrow{e}\node{\Cl_K}\arrow{e}\arrow{s}\node{M_{(\varepsilon)}}\arrow{e}
\arrow{s}\node{I_G}\arrow{e}\node{0}\\
\node[2]{0}\node{0}
\end{diagram}
\end{equation*}
where $(\overline{\varepsilon})$ is the group extension
\[
1\longrightarrow C_K
\longrightarrow \overline{\mathcal{G}}
\longrightarrow G
\longrightarrow 1
\]
associated to the fundamental class $c_{K/k}\in H^2(G,C_K)$,
and 
$(\varepsilon)$ is the group extension
\[
1\longrightarrow \Cl_K
\longrightarrow \mathcal{G}
\longrightarrow G
\longrightarrow 1 
\]
associated to the image of $c_{K/k}$ under the natural map
$H^2(G,C_K)\longrightarrow H^2(G,\Cl_K)$.
It follows from Shafarevich's theorem \cite[Chapter 15, Theorem 6]{AT} that
the group extensions $(\varepsilon)$ and $(\varepsilon_{K/k})$
of $G$ are isomorphic via a morphism inducing
the Artin map $\Cl_K\simeq\Gal(H_K/K)$,
hence  
$M_{(\varepsilon)}\simeq M_{(\varepsilon_{K/k})}$
as $\Z[G]$-modules.
The $\Z[G]$-module $U_K$ is cohomologically trivial since
$K/k$ is unramified, and
we know $M_{(\overline{\varepsilon})}$
is also cohomologically trivial
(see, for example, \cite[Theorem (3.1.4)]{NSW}).
Therefore we obtain the $\Z_p[G]$-module isomorphisms
\begin{align*}
\varphi_H:\hat{H}^i(H,E_K)&\simeq \hat{H}^{i-1}(H, U_K/E_K)
\simeq \hat{H}^{i-1}(H,N)\\
&\simeq \hat{H}^{i-2}(H,M_{(\varepsilon)})\simeq
\hat{H}^{i-2}(H,M_{(\varepsilon_{K/k})}),
\end{align*}
for the subgroups $H$ of $G$ from the above diagram.
Furthermore we can see that the maps $\varphi_H$'s are compatible with the 
restriction maps and the co-restriction maps.
\hfill$\Box$

\

\begin{rem}
In the case where $G=\Gal(K/k)$ is a $p$-group, 
we define the group
extension
$(\varepsilon_{K/k})_p$ to be
\[
(\varepsilon_{K/k})_p\ \
1\longrightarrow\Gal(H_K(p)/K)\longrightarrow
\Gal(H_K(p)/k)\longrightarrow\Gal(K/k)\longrightarrow 1,
\]
where $H_K(p)$ is the maximal unramified abelian $p$-extension field of $K$.
Under this situation, we obtain the natural $\Z_p[G]$-module isomorphisms
\[
\hat{H}^i(H,M_{(\varepsilon_{K/k})})\simeq \hat{H}^i(H,M_{(\varepsilon_{K/k})_p}).
\]
for the subgroups $H$ of $G$, which are compatible with the restriction
and co-restriction maps.
\end{rem}
In the rest of this section, we will show that there exists
an unramified $G$-extension $K/k$ such that
the cohomology groups of $\mathcal{E}_K$
are isomorphic to those of a given $\Z_p[G]$-module $C$
stated in Theorem 1.
We need some lemmas:
\begin{lem}
Let $G$ be a finite $p$-group and 
\[
0\longrightarrow A\longrightarrow B\overset{\pi}{\longrightarrow} I_{\Z_p[G]}\longrightarrow 0
\]
be a exact sequence of $\Z_p[G]$-modules with finite $A$,
where $I_{\Z_p[G]}$ stands for the augumentation ideal of $\Z_p[G]$.
Then there exists an unramified $G$-extension of number fields
$K/k$ such that
\[
B\simeq M_{(\varepsilon_{K/k})_p}\otimes_{\Z}\Z_p
\]
as $\Z_p[G]$-modules.
Furthermore we can choose $K/k$ so that
$[k:\Q]$ is arbitrarily large, and that $\mu_p\not\subseteq K$
if $p\ne 2$,
$\mu_p$ being the $p$-th roots of unity.
\end{lem}
{\it Proof.}\ \ \ 
Let $s$ be a fixed section of $\pi$ as $\Z_p$-modules and
define the map $f:G\times G\longrightarrow A$
by
\[
f(\sigma,\tau):=
\sigma s(\tau-1)-s(\sigma\tau-1)+s(\sigma-1)\in A\ \ \ (\sigma,\tau\in G).
\]
Then we find that $f$ is a 2-cocycle.
Furtheermore if
\[
(\varepsilon)\ \ 1\longrightarrow A\longrightarrow\mathcal{G}\longrightarrow G\longrightarrow 1
\]
is the group extension of $G$ by $A$ associated to the cohomology class
$[f]\in H^2(G,A)$, 
then we see that 
\[
M_{(\varepsilon)}\otimes_{\Z}\Z_p\simeq B
\]
as $\Z_p[G]$-modules.

Since $\mathcal{G}$ is a $p$-group, 
it follows from \cite{HMR}, \cite{Oz} that
there exists a number field $k$ such that we have an isomorphism
\[
\psi:\mathcal{G}\simeq\Gal(L_p(k)/k),
\]
where $L_p(k)$ is the maximal unramified $p$-extension of $k$. 
Here we may choose $K/k$ so that $[k:\Q]$ is arbitrarily large, and that $\mu_p\not\subseteq K$ if $p\ne 2$ (See \cite{HMR}, \cite{Oz}).
Put $K:=L_p(k)^{\psi(A)}$. Then $K/k$ is an unramified $G$-extension
and the group extension
\[
(\varepsilon_{K/k})_p\ \ \ 1\longrightarrow\Gal(L_p(k)/K)\longrightarrow\Gal(L_p(k)/k)
\longrightarrow\Gal(K/k)\longrightarrow 1
\]
is isomorphic to $(\varepsilon)$.
Hence we have
\[
M_{(\varepsilon_{K/k})_p}\simeq
M_{(\varepsilon)}
\]
as $\Z[G]$-modules.
Thus we conclude that 
\[
M_{(\varepsilon_{K/k})_p}\otimes_\Z\Z_p\simeq
M_{(\varepsilon)}\otimes_\Z\Z_p\simeq B
\]
as $\Z_p[G]$-modules.
\hfill$\Box$

\

\begin{lem}
Let $G$ be a finite cyclic $p$-group and 
$C$ a finitely generated 
$\Z_p[G]$-module such that
\[
C\otimes_{\Z_p}\Q_p\simeq\Q_p[G]^{\oplus r}\oplus I_{\Q_p[G]}.
\]
Then there exists an exact sequence of $\Z_p[G]$-modules
\[
0\longrightarrow A\longrightarrow B\overset{\pi}{\longrightarrow} I_{\Z_p[G]}\longrightarrow 0
\]
with finite $A$ and $\Z_p[G]$-module isomorphisms
\[
\psi_H:\hat{H}^i(H,C)\simeq\hat{H}^{i-2}(H,B)
\]
for the subgroups $H$ of $G$ and $i\in\Z$
such that $\psi_H$'s are compatible with the restriction
and the co-restriction maps.
\end{lem}
{\it Proof.}\ \ \ 
It follows from the assumption on $C$ that
$C$ has a finite index $\Z_p[G]$-submodule 
\[
F\oplus I
\]
such that $F\simeq\Z_p[G]^{\oplus r}$ and $I\simeq I_{\Z_p[G]}$
as $\Z_p[G]$-modules.
Put $C_1:=C/F$.
Then we see that $C_1$ has a finite index $\Z_p[G]$-submodule $I_1$
isomorphic to $I_{\Z_p[G]}$ as $\Z_p[G]$-modules.
Let $x_0\in I_1$ be a $\Z_p[G]$-generator of $I_1$
corresponding to $\sigma-1\in I_{\Z_p[G]}$ under
$I_1\simeq I_{\Z_p[G]}$, $\sigma\in G$ being a fixed generator of
the cyclic group $G$.
Define $C_2:=C_1\oplus\Z_p[G]$
and the $\Z_p[G]$-homomorphism
\[
\iota:\Z_p[G]\longrightarrow C_2=C_1\oplus\Z_p[G]
\]
by
\[
\iota(\alpha):=(\alpha x_0,\alpha N_G),
\]
where $N_G:=\sum_{\tau\in G}\tau\in\Z_p[G]$.
Then we see that $\iota$ is injective since
$\alpha x_0=0\Longleftrightarrow\alpha\in \Z_p N_G$
and $\alpha N_G=0\Longleftrightarrow\alpha\in I_{\Z_p[G]}$ hold.
Put $B:=C_2/\iota(\Z_p[G])$.
Then we get the $\Z_p[G]$-isomorphisms
\begin{equation}\label{isom1}
\hat{H}^i(H,C)\simeq\hat{H}^i(H,C_1)\simeq\hat{H}^i(H,C_2)\simeq
\hat{H}^i(H,B)
\end{equation}
for the subgroups $H$ of $G$, which are compatible with
the restriction and the co-restriction maps.
Furthermore, since $G$ is a cyclic group,
there exist $\Z_p[G]$-module isomorphisms
\begin{equation}\label{isom2}
\hat{H}^i(H,B)\simeq\hat{H}^{i-2}(H,B)
\end{equation}
for the subgroups $H$ of $G$ 
such that they are compatible with the restriction and the co-restriction maps.

We define the map
\[
\pi:B=C_2/\iota(\Z_p[G])\longrightarrow I_{\Z_p[G]}
\]
by $\pi((x,y)+\iota(\Z_p[G])):=(\sigma-1)y$
for $(x,y)\in C_2=C_1\oplus\Z_p[G]$.
Then $\pi$ is a surjective $\Z_p[G]$-homomorphism, and 
we get the exact sequence of $\Z_p[G]$-modules
\[
0\longrightarrow A\longrightarrow B\overset{\pi}{\longrightarrow} I_{\Z_p[G]}\longrightarrow 0
\]
with $A:=\ker\pi$.
Here $A$ is finite since 
\[
\mathrm{rank}_{\Z_p}B
=\rank_{\Z_p}C_1=\rank_{\Z_p}I_1=\rank_{\Z_p}{I_{\Z_p[G]}}.
\]
Thus we have completed the proof by \eqref{isom1} and \eqref{isom2}.
\hfill$\Box$

\

Combining Propositions 1 and Lemmas 1 and 2, we obtain the following:
\begin{prop}
Let $p$ be an odd  prime number and $G$ a finite cyclic $p$-group.
Assume that $C$ is any given finitely generated $\Z_p[G]$-module
with $\mathrm{Tor}_{\Z_p}(C)=0$ and
\[
C\otimes_{\Z_p}\Q_p\simeq \Q_p[G]^{\oplus r}\oplus 
I_{\Q_p[G]}.
\]
Then there exists an unramified $G$-extension $K/k$ of number fields
equipped with $\Z_p[G]$-isomorphisms
\[
\hat{H}^i(H,\mathcal{E}_K)\simeq\hat{H}^i(H,C)\ \  (i\in\Z)
\]
for the subgroups $H$ of $G$
which are compatible with the restriction and the co-restriction maps.
We may choose $K/k$ so that $[k:\Q]$ is arbitrarily large.
\end{prop}

{\bf Proof.}\ \ \ 
It follows from Lemma 2
that there exists an exact sequence of $\Z_p[G]$-modules
\[
0\longrightarrow A\longrightarrow B\overset{\pi}{\longrightarrow} I_{\Z_p[G]}\longrightarrow 0
\]
with finite $A$ and $\Z_p[G]$-module isomorphisms
\[
\psi_H:\hat{H}^i(H,C)\simeq\hat{H}^{i-2}(H,B)
\]
for the subgroups $H$ of $G$ and $i\in\Z$
such that $\psi_H$'s are compatible with the restriction
and the co-restriction maps.

By using Lemma 1, we have an unramified $G$-extension $K/k$
of number fields such that
\[
B\simeq M_{(\varepsilon_{K/k})_p}\otimes_{\Z}\Z_p
\]
as $\Z_p[G]$-modules via a fixed identification of $\Gal(K/k)$ with $G$.
Here We can choose $K/k$ so that $[k:\Q]$ is arbitrarily large, and that $\mu_p\not\subseteq K$.

Therefore, by using Proposition 1 and Remark 1,
we obtain the $\Z_p[G]$-module isomorphisms
\begin{align*}
\hat{H}^i(H,C)&\simeq\hat{H}^{i-2}(H,B)
\simeq\hat{H}^{i-2}(H,M_{(\varepsilon_{K/k})_p}\otimes_{\Z}\Z_p)\\
&\simeq\hat{H}^{i-2}(H,M_{(\varepsilon_{K/k})_p})
\simeq\hat{H}^{i-2}(H,M_{(\varepsilon_{K/k})})\\
&\simeq\hat{H}^i(H,E_K)
\simeq\hat{H}^i(H,\mathcal{E}_K)
\end{align*}
for the subgroups $H$ of $G$, 
which are compatible with the restriction and the co-restriction
maps.
\hfill$\Box$
\section{Yakovrev's theory}
Let $p$ be a prime number and $G$ a cyclic $p$-group.
In this section we introduce Yakovrev's theory \cite{Ya}
for the determination of $\Z_p$-torsion free $\Z_p[G]$-modules
in terms of a certain
diagram consisting of 1-dimensional cohomology groups
and the restriction and the co-restriction maps.

Let $G$ be a cyclic group of order
$p^n\ (n\ge 1)$.
We write $G_i$ for the subgroup of $G$ of order $p^i\ (0\le i\le n)$.
We define $\frak{M}_G$ to be the category of the diagrams
\[
(A_\bullet,\alpha_\bullet,\beta_{\bullet}):\ A_1\overset{\alpha_1}{\underset{\beta_1}{\rightleftarrows}} A_2
\overset{\alpha_2}{\underset{\beta_2}{\rightleftarrows}}\cdots
\overset{\alpha_{n-2}}{\underset{\beta_{n-2}}{\rightleftarrows}}
A_{n-1}\overset{\alpha_{n-1}}{\underset{\beta_{n-1}}
{\rightleftarrows}}A_n
\]
in which each $A_i$ is a finite $(\Z/p^i)[G/G_i]$-module, 
and each $\alpha_i$ and $\beta_i$ is a $\Z_p[G]$-module homomorphism
such that $\alpha_i\circ\beta_i$ and $\beta_i\circ\alpha_i$ are respectively induced by multiplication by $p$ and by the
action of $\sum_{g\in G/G_i}g$. 

A morphism $(A_\bullet, \alpha_\bullet,\beta_\bullet)\longrightarrow
(A'_\bullet, \alpha'_\bullet,\beta'_\bullet)$
in $\frak{M}_G$ is defined to be
a collection of $\Z_p[G]$-module homomorphisms
$\{\gamma_i : A_i\longrightarrow A'_i\,|\,1\le i\le n\}$
which is compatible with the respective maps $\alpha_\bullet,\, \beta_\bullet,\, 
\alpha'_\bullet,\,\beta'_\bullet$.
In particular, such a morphism is an isomorphism if and only if each 
map $\gamma_i$ is bijective.

For any finitely generated $\Z_p[G]$-module $M$,
we define $\Delta(M)\in\mathfrak{M}_G$
to be the diagram
\[
H^1(G_1,M)\overset{\alpha_1}{\underset{\beta_1}{\rightleftarrows}}
H^1(G_2,M)\overset{\alpha_2}{\underset{\beta_2}{\rightleftarrows}}\cdots
\overset{\alpha_{n-2}}{\underset{\beta_{n-2}}{\rightleftarrows}}
H^1(G_{n-1},M)\overset{\alpha_{n-1}}{\underset{\beta_{n-1}}
{\rightleftarrows}}H^1(G_n,M),
\]
where $\alpha_i$'s and $\beta_i$'s are the co-restriction maps and 
the restriction maps, respectively.

The following plays a crucial role in the proof of Theorem 1:
\begin{prop}[{\cite[Theorem 4]{Ya}}]
Let $M$ and $M'$ be finitely generated $\Z_p[G]$-modules
with $\mathrm{Tor}_{\Z_p}(M)=\mathrm{Tor}_{\Z_p}(M')=0$.
If $\Delta(M)\simeq\Delta(M')$ in $\mathfrak{M}_G$ and $H^2(G,M)\simeq H^2(G,M')$ holds,
then 
\[
M\oplus\Z_p[G]^a\simeq M'\oplus\Z_p[G]^b
\]
as $\Z_p[G]$-modules holds for certain integers $a,b\ge 0$.
\end{prop}
\section{Proof of Theorem 1}
Now we will give a proof of Theorem 1.

\

{\bf Proof of Theorem 1}\ \ \ 
Let $p$ be an odd prime number and $G$ a finite cyclic group
of order $p^n$.
Assume that $C$ is any given finitely generated $\Z_p[G]$-module
with $\mathrm{Tor}_{\Z_p}(C)=0$ and
\[
C\otimes_{\Z_p}\Q_p\simeq \Q_p[G]^{\oplus r}\oplus 
I_{\Q_p[G]}.
\]
Then it follows from Proposition 2
that there exists an unramified $G$-extension $K/k$ of number fields
equipped with $\Z_p[G]$-module isomorphisms
\[
\hat{H}^i(H,\mathcal{E}_K)\simeq\hat{H}^i(H,C)\ \  (i\in\Z)
\]
for the subgroups $H$ of $G$
which are compatible with the restriction and the co-restriction maps.
This implies that 
\[
\Delta(\mathcal{E}_K)\simeq\Delta(C)
\]
in the category $\frak{M}_G$
and
\[
H^2(G,\mathcal{E}_K)\simeq H^2(G,C)
\]
hold.

Hence we conclude by Proposition 3 that
\[
C\oplus\Z_p[G]^{\oplus a}\simeq\mathcal{E}_K\oplus\Z_p[G]^{\oplus b}
\]
as $\Z_p[G]$-modules for certain integers $a,\, b\ge 0$.
If we choose $K/k$ 
with large $[k:\Q]$ so that $\rank_{\Z_p}\mathcal{E}_K\ge\rank_{\Z_p}C$,
then we have $b\le a$.
Therefore we see that
\[
\mathcal{E}_K\simeq C\oplus\Z_p[G]^{\oplus(a-b)}
\]
as $\Z_p[G]$-modules holds by the Krull-Schmidt theorem for the finitely generated
$\Z_p[G]$-modules.
\hfill$\Box$
\section{Examples of explicit $\Z_p[G]$-module structures of $\mathcal{E}_K$}
Let $p$ be an odd prime number and $K/k$ a unramified cyclic 
$p$-extension with $G=\Gal(K/k)$ and $\mu_p\not\subseteq K$. 

From Proposition 1, Remark 1, and Proposition 3,
we conclude that the $\Z_p[G]$-module structure of $\mathcal{E}_K$
modulo free direct factors
is completely determined by the group extension
\[
(\varepsilon_{K/k})_p\ \ \ 1\longrightarrow
\Gal(H_K(p)/K)\longrightarrow\Gal(H_K(p)/k)
\longrightarrow G\longrightarrow 1.
\]

We will close this paper by giving
examples of explicit $\Z_p[G]$-module structures
of $\mathcal{E}_K$ for a certain class of $K/k$
by using the above fact.

The first example is the special case $S=\emptyset$
of \cite[Theorem 4.4]{BLM}:
\begin{prop}
Let $p$ be an odd prime number and $G$ a cyclic $p$-group.
Assume that $K/k$ is an unramified $G$-extension such
that $\Gal(H_K(p)/k)$ is cyclic and $\mu_p\not\subseteq k$.
Then we have
\[
\mathcal{E}_K\simeq I_{\Z_p[G]}\oplus\Z_p[G]^{\oplus
\mathrm{rank}_{\Z_p}\mathcal{E}_k}
\]
as $\Z_p[G]$-modules.
\end{prop}
{\bf Proof}\ \ \ 
Put $L:=H_K(p)$ and $\mathcal{G}:=\Gal(L/k)$.
Then $\mathcal{G}$ is a cyclic $p$-group by the assumption on $K/k$.
We note that $H_L(p)=L$ by the Burnside basis theorem.
We will apply Proposition 1 to the unramified cyclic $p$-extension $L/k$.
We note that
\[
M_{(\varepsilon_{L/k})_p}\otimes_{\Z}\Z_p\simeq I_{\Z_p[\mathcal{G}]}
\]
as $\Z_p[\mathcal{G}]$-modules since $\Gal(H_L(p)/L)=1$.
Hence we have the $\Z_p[\mathcal{G}]$-isomorphisms
\begin{align*}
\hat{H}^i(\mathcal{H},\mathcal{E}_L)&\simeq
\hat{H}^i(\mathcal{H},E_L)\simeq
\hat{H}^{i-2}(\mathcal{H},M_{(\varepsilon_{L/k})_p}\otimes_{\Z}\Z_p)\\
&\simeq
\hat{H}^{i-2}(\mathcal{H},I_{\Z_p[\mathcal{G}]})\simeq
\hat{H}^i(\mathcal{H},I_{\Z_p[\mathcal{G}]})
\end{align*}
(the last isomorphism exists since $\mathcal{G}$ is cyclic)
for the subgroups $\mathcal{H}$ of $\mathcal{G}$ which
are compatible with the restriction and the co-restriction maps.
This implies 
\[
\Delta(\mathcal{E}_L)\simeq\Delta(I_{\Z_p[\mathcal{G}]})
\]
in the category $\mathfrak{M}_{\mathcal{G}}$ and
\[
H^2(\mathcal{G},\mathcal{E}_L)
\simeq H^2(\mathcal{G},I_{\Z_p[\mathcal{G}]}).
\]
Hence we conclude that
\[
\mathcal{E}_L\simeq I_{\Z_p[\mathcal{G}]}
\oplus\Z_p[\mathcal{G}]^{\oplus\rank_{\Z_p}\mathcal{E}_k}
\]
as $\Z_p[\mathcal{G}]$-modules by Proposition 3 and Dirichlet's unit theorem.
Then we derive the $\Z_p[G]\simeq\Z_p[\mathcal{G}/\Gal(L/K)]$-module
isomorphism
\begin{align*}
\mathcal{E}_K=\mathcal{E}_L^{\Gal(L/K)}
&\simeq I_{\Z_p[\mathcal{G}]}^{\Gal(L/K)}\oplus
\left(\Z_p[\mathcal{G}]^{\Gal(L/K))}\right)
^{\oplus\rank_{\Z_p}\mathcal{E}_k}\\
&\simeq I_{\Z_p[G]}\oplus\Z_p[G]^{\mathrm{rank}_{\Z_p}\mathcal{E}_k}
\end{align*}
from the above isomorphism of $\Z_p[\mathcal{G}]$-modules.
\hfill$\Box$

\

To give the second example, we introduce a certain $\Z_p[G]$-module:
Let $G$ be a cyclic group of order $p^n\ (n\ge 1)$.
We define
\[
J_e:=p^e\Z_p[G]+\Z_p N_G\subseteq\Z_p[G]
\]
for $e\ge 0$. We note that $J_0=\Z_p[G]$.

We first give the following lemma for the $\Z_p[G]$-modules $J_e$:
\begin{lem}
(1)\ \ \ We have isomorphisms of $\Z_p[G]$-modules
\[
\hat{H}^i(H,J_e)\simeq\hat{H}^i(H,\Z/p^e)
\]
for the subgroups $H$ of $G$ which are compatible
with the restriction and co-restriction maps.

Furthermore, we have
\[
\hat{H}^{i}(G,\Z/p^e)\simeq\Z/{p^{\min\{e,n\}}}
\]
as $\Z_p[G]$-modules. 

(2)\ \ \ The $\Z_p[G]$-modules $J_e$'s are not isomorphic to each other for $0\le e\le n$.
Also
\[
J_e\simeq\Z_p\oplus I_{\Z_p[G]}
\]
as $\Z_p[G]$-modules holds for any $e\ge n$.
\end{lem}
{\bf Proof.}\ \ \ 
(1)\ \ \ We define the surjective $\Z_p[G]$-module homomorphism
\[
f:\Z_p[G]\longrightarrow\Z/p^e
\]
by
\[
f(\alpha):=\mathrm{aug}(\alpha)+p^e\Z,
\]
$\mathrm{aug}$ being the augumentation map on $\Z_p[G]$.
Then we see that
\[
\mathrm{ker}\,f=p^e\Z_p[G]+I_{\Z_p[G]}.
\]
Furthermore we define the surjective $\Z_p[G]$-module homomorphism
\[
g:\Z_p[G]^{\oplus 2}\longrightarrow \mathrm{ker}\,f
\]
by 
\[
g((\alpha,\beta)):=p^e\alpha+(\sigma-1)\beta,
\]
$\sigma$ being a fixed generator of $G$.
Then we find that
\[
\mathrm{ker}\,g=
\{((\sigma-1)\alpha,-p^e\alpha+aN_G)\,|\,\alpha\in\Z_p[G], a\in\Z_p\}.
\]
Here we have the $\Z_p[G]$-module isomorphism
\[
\mathrm{ker}g\simeq J_e
\]
defined by
\[
((\sigma-1)\alpha, -p^e\alpha+aN_G)\mapsto -p^e\alpha+aN_G.
\]
Thus we get the exact sequence
\[
0\longrightarrow J_e{\longrightarrow}\Z_p[G]^{\oplus 2}
\overset{g}{\longrightarrow}\Z_p[G]\overset{f}{\longrightarrow}\Z/p^e
\longrightarrow 0
\]
of $\Z_p[G]$-modules.
This induces the $\Z_p[G]$-module isomorphisms
\[
\hat{H}^i(H,J_e)\simeq\hat{H}^{i-2}(H,\Z/p^e)\simeq\hat{H}^i(H,\Z/p^e)
\]
for the subgroups $H$ of $G$ which are
compatible with the restriction and co-restriction
maps since $G$ is a cyclic group.

We can get $\hat{H}^i(G,\Z/p^e)\simeq\Z/p^{\min\{e,n\}}$ by the exact sequence
\[
0\longrightarrow\Z\overset{p^e}{\longrightarrow}\Z
\longrightarrow\Z/p^e\longrightarrow 0
\]
and the fact that
\begin{equation*}
\hat{H}^i(G,\Z)\simeq
\begin{cases}
0\ \ (\mbox{if $i$ is odd} ),\\
\Z/p^n\ \ (\mbox{if $i$ is even}).
\end{cases}
\end{equation*}

\

(2)\ \ \ The first assertion is obvious from
\[
\hat{H}^i(G,J_e)\simeq\Z/p^e
\]
for $0\le e\le n$, which follows from 
(1) for $H=G$.

We will give an explicit $\Z_p[G]$-module isomorhism 
\[
h: \Z_p\oplus I_{\Z_p[G]}\simeq J_e.
\]
for $e\ge n$.
Define 
\[
h(a,(\sigma-1)\beta):=p^e\beta-p^{e-n}\mathrm{aug}(\beta)N_G+aN_G,
\]
for $a\in\Z_p,\beta\in\Z_p[G]$,
where $\sigma$ is a fixed generator of $G$.
We can see that the above $h$ is well-defined $\Z_p[G]$-homomorphism
since
$(\sigma-1)\beta=(\sigma-1)\beta'$ for $\beta,\beta'\in\Z_p[G]$
implies $\beta-\beta'\in\Z_pN_G$.
Surjectivity of $h$ is obvious.
Furthermore we find that $h$ is injective
because $p^e\beta-p^{e-n}\mathrm{aug}(\beta)N_G+aN_G=0$
implies $p^e\beta\in p^e\Z_p[G]\cap\Z_pN_G=p^e\Z_pN_G$,
which in turn means $\beta\in\Z_pN_G$ and $aN_G=-(p^e\beta-p^{e-n}
\mathrm{aug}(\beta)N_G)=0$.
Thus we conclude that $h$ is a $\Z_p[G]$-module isomorphism.
\hfill$\Box$

\

Now we will give the second example:
\begin{prop}
Let $p$ be an odd prime number and $G$ a cyclic $p$-group of order
$p^n$.
Assume that $K/k$ is an unramified $G$-extension such
that $\mathcal{G}:=\Gal(H_p(K)/k)$ is abelian and the group extension
\[
(\varepsilon_{K/k})_p\ \ \ 1\longrightarrow
\Gal(H_K(p)/K)\longrightarrow\mathcal{G}
\longrightarrow G\longrightarrow 1.
\]
splits, and that
\[
\Gal(H_K(p)/K)\simeq\bigoplus_{j=1}^r\Z/p^{e_j}
\]
as abelian groups for $e_j\ge 1$.
Then we have 
\[
\mathcal{E}_K\simeq\bigoplus_{j=1}^r J_{\min\{e_j,n\}}
\oplus I_{\Z_p[G]}\oplus\Z_p[G]^{\oplus(\rank_{\Z_p}\mathcal{E}_k-r)}
\]
as $\Z_p[G]$-modules.
\end{prop}

\

{\bf Proof.}\ \ \ 
It follows from the assumption on $K/k$ that
\[
M_{(\varepsilon_{K/k})_p}\otimes_{\Z}\Z_p\simeq
\bigoplus_{j=1}^r\Z/p^{e_j}\oplus I_{\Z_p[G]}
\]
as $\Z_p[G]$-modules since
the action of $G$ on $\Gal(H_K(p)/K)$
and the extension class of $(\varepsilon_{K/k})_p$ are trivial.
Then, from Proposition 1 and Lemma 3, we have
$\Z_p[G]$-module isomorphisms
\begin{align*}
\hat{H}^i(H,\mathcal{E}_K)&\simeq
\hat{H}^{i-2}(H,M_{(\varepsilon_{K/k})_p}\otimes_{\Z}\Z_p)\\
&\simeq
\bigoplus_{j=1}^r\hat{H}^{i-2}(H,\Z/p^{e_j})
\oplus\hat{H}^{i-2}(H,I_{\Z_p[G]})\\
&\simeq
\bigoplus_{j=1}^r\hat{H}^{i-2}(H,J_{e_j})
\oplus\hat{H}^{i-2}(H,I_{\Z_p[G]})\\
&\simeq\hat{H}^i(H,\bigoplus_{j=1}^r J_{\min\{e_j,n\}}
\oplus I_{\Z_p[G]})
\end{align*}
for the subgroups $H$ of $G$
which are compatible with the restriction and the co-restriction maps. 
Hence 
\[
\Delta(\mathcal{E}_K)\simeq\Delta(\bigoplus_{j=1}^r J_{\min\{e_j,n\}}
\oplus I_{\Z_p[G]})
\]
in $\mathfrak{M}_G$ and
\[
H^2(G,\mathcal{E}_K)\simeq H^2(G,\bigoplus_{j=1}^r J_{\min\{e_j,n\}}
\oplus I_{\Z_p[G]})
\]
hold.
Therefore we conclude that
\[
\mathcal{E}_K\simeq
\bigoplus_{j=1}^r J_{\min\{e_j,n\}}
\oplus I_{\Z_p[G]}\oplus\Z_p[G]^{\oplus(\rank_{\Z_p}\mathcal{E}_k-r)}
\]
as $\Z_p[G]$-modules holds by using Proposition 3 and Dirichlet's unit theorem.
\hfill$\Box$

Manabu Ozaki,\par\noindent
Department of Mathematics,\par\noindent
School of Fundamental Science and Engineering,\par\noindent
Waseda University,\par\noindent
Ohkubo 3-4-1, Shinjuku-ku, Tokyo, 169-8555, Japan\par\noindent
e-mail:\ \verb+ozaki@waseda.jp+
\end{document}